\documentstyle[amscd,amssymb,verbatim,epsf]{amsart}

\begin{document}
\theoremstyle{plain}
\newtheorem{Thm}{Theorem}
\newtheorem{Cor}{Corollary}
\newtheorem{Con}{Conjecture}
\newtheorem{Main}{Main Theorem}
\newtheorem{Lem}{Lemma}
\newtheorem{Prop}{Proposition}

\theoremstyle{definition}
\newtheorem{Def}{Definition}
\newtheorem{Note}{Note}
\newtheorem{Ex}{Example}

\theoremstyle{remark}
\newtheorem{notation}{Notation}
\renewcommand{\thenotation}{}

\errorcontextlines=0
\numberwithin{equation}{section}
\renewcommand{\rm}{\normalshape}%

\title
   {The Geometry of Focal Sets}

\author{Brendan Guilfoyle}
\address{Brendan Guilfoyle\\
          Department of Mathematics and Computing \\
          Institute of Technology, Tralee \\
          Clash \\
          Tralee  \\
          Co. Kerry \\
          Ireland.}
\email{brendan.guilfoyle@@ittralee.ie}

\author{Wilhelm Klingenberg}
\address{Wilhelm Klingenberg\\
 Department of Mathematical Sciences\\
 University of Durham\\
 Durham DH1 3LE\\
 United Kingdom.}
\email{wilhelm.klingenberg@@durham.ac.uk }

\date{October 23rd, 2005}

\begin{abstract}
The space ${\Bbb{L}}$ of oriented lines, or rays, in ${\Bbb{R}}^3$ is a 4-dimensional space with an abundance of natural 
geometric structure. In particular, it boasts a neutral K\"ahler metric which is closely related to the Euclidean metric on 
${\Bbb{R}}^3$. In this paper we explore the relationship between the focal set of a line congruence (or 2-parameter family 
of oriented lines in ${\Bbb{R}}^3$) and the geometry induced on the associated surface in ${\Bbb{L}}$. The physical context
of such sets is geometric optics in a homogeneous isotropic medium, and so, to illustrate the method, we 
compute the focal set
of the $k^{\mbox{th}}$ reflection of a point source off the inside of a cylinder. The focal sets, which we 
explicitly parameterize, exhibit unexpected symmetries, and are found to fit well with observable phenomena.

\end{abstract}

\maketitle

The space of oriented affine lines in ${\Bbb{R}}^3$, which we denote ${\Bbb{L}}$, has an
abundance of natural geometric structure. The purpose of this paper is to continue recent work 
\cite{gak3} relating this structure
to geometric optics in a homogeneous isotropic medium: the theory of light propagation under the assumption that the light travels 
along straight lines in ${\Bbb{R}}^3$. The fundamental objects of study are 2-parameter families of oriented lines, or 
{\it line congruences}, which we view as surfaces in ${\Bbb{L}}$. Thus we are lead to consider the
geometry of immersed surfaces $\Sigma\subset{\Bbb{L}}$.

In the first instance, since ${\Bbb{L}}$ can be identified with the tangent space to the 2-sphere, there is the 
natural bundle map $\pi:{\Bbb{L}}\rightarrow S^2$. If 
$\pi|_{\Sigma}:\Sigma\rightarrow S^2$ is not an immersion, we say that $\Sigma$ is {\it flat}. Otherwise, 
$\Sigma$ can be described, at least locally, by sections of the canonical bundle.

On the other hand, there is a natural symplectic structure $\Omega$ on ${\Bbb{L}}$, and $\Sigma\subset{\Bbb{L}}$ is lagrangian 
with respect to this symplectic structure iff the line congruence admits a family of orthogonal surfaces in 
${\Bbb{R}}^3$. In geometric optics such surfaces are the wavefronts of the propagating light. 
As a wavefront evolves along the line congruence, if there is any focusing, the surface becomes singular. The points at 
which this occurs are referred to as {\it focal points}.

In addition, ${\Bbb{L}}$ admits a natural complex structure ${\Bbb{J}}$, which, together with the symplectic structure, 
forms a
natural K\"ahler structure \cite{gak4}.  The metric ${\Bbb{G}}$ is of 
signature ($++--$) and therefore the metric induced on a surface $\Sigma$ may be riemannian, lorentz or degenerate.  

The aim of this paper is two-fold: to relate the geometry induced on $\Sigma\subset{\Bbb{L}}$ by ${\Bbb{G}}$ with the set of focal 
points of the line congruence in ${\Bbb{R}}^3$,
and to demonstrate the effectiveness of the geometric formalism in computing focal sets explicitly. Moreover, we hope
to demonstrate that the results of such computations can have visible physical significance.

For the first aim we prove:

\vspace{0.1in}

\noindent {\bf Main Theorem 1.}

{\it
Let $\Sigma$ be an immersed surface in ${\Bbb{L}}$. If $\Sigma$ is flat, there is exactly one focal point on each line of
the congruence. If $\Sigma$ is not flat, then there is none, one or two focal points on each line iff the
metric induced on $\Sigma$ by ${\Bbb{G}}$ is riemannian, degenerate or lorentz (respectively).
}

\vspace{0.1in}
In addition, we relate the distance between pairs of focal surfaces and the angle between their normals to 
geometric quantities on $\Sigma\subset{\Bbb{L}}$.

Mathematically, the focal set of a generic line congruence is well 
understood \cite{agv} \cite{bgg2} \cite{ist}. Special examples of focal sets, also referred to as caustics 
by some authors, have been studied for many decades \cite{bgg1} \cite{caly} \cite{glae1} \cite{hold}. 
For the second aim, we compute explicitly the focal set of the line congruence formed by the multiple reflection of a 
point source off the inside of a cylinder.  The first such reflected focal set is often referred to as the coffeecup 
caustic, since its cross-section is commonly observed on the top of a cup of coffee in the presence of a bright
light. We prove:

\vspace{0.1in}

\noindent {\bf Main Theorem 2.}

{\it
Consider the k$^{th}$ reflection of a point source at ($-l,0,0$) off the inside of a cylinder lying along the x$^3$-axis with
radius $a$.
The focal set of the reflected line congruence consists of 
a surface:
\[
z=(-1)^{k+1}l\frac{[l\sin v\;i+(a^2-l^2\sin^2v)^{\frac{1}{2}}]^{2k}
     [2kl\cos v\sin v\;e^{-iv}+(a^2-l^2\sin^2v)^{\frac{1}{2}}]}
     {a^{2k}[2kl\cos v+(a^2-l^2\sin^2v)^{\frac{1}{2}}]}
\]
\[
x^3=\frac{k(1-u^2)[a^2-l^2-2l^2\sin^2v+2kl\cos v(a^2-l^2\sin^2v)^{\frac{1}{2}}]}{u[2kl\cos v+(a^2-l^2\sin^2v)^{\frac{1}{2}}]}
\]
and a curve in the $x^1x^2$-plane: 
\[
z=(-1)^{k+1}ka^{-2k}[l\sin v\;i+(a^2-l^2\sin^2v)^{\frac{1}{2}}]^{2k}
  [l+2ke^{iv}(a^2-l^2\sin^2v)^{\frac{1}{2}}]
\]
where $z=x^1+ix^2$, $u\in{\Bbb{R}}$ and $v$ is in the domain
\[
0\le v\le\pi \qquad\qquad\mbox{for} \qquad l\le a 
\]
and
\[-\sin^{-1}(a/l)\le v\le \sin^{-1}(a/l)
\qquad\qquad\mbox{for}\qquad l> a.
\]
}

\vspace{0.1in}

Thus the focal set consists of two sets: a translation invariant surface and a curve lying in the plane
perpendicular to the symmetry axis and containing the source point. The commonly observed coffeecup caustic is the 
level sets of the former with $k=1$, while the latter lies outside the cup and so is not observed. With the aid of a 
polished brass cylinder, the first author, in collaboration with Grace Weir, has photographed multiple reflections 
of up to 4$^{\mbox{th}}$ order. The results agree closely with the result of Main Theorem 2.

This paper is organised as follows. The next section contains a review of the complex geometric structure on 
the space of oriented lines in
${\Bbb{R}}^3$, as developed in \cite{gak4}. The second section describes the first order invariants of 
an arbitrary line congruence, while the next section contains an exposition on the construction of surfaces
using line congruences. This is required in order to have a description of the reflective surfaces. To illustrate
the method we compute the normal line congruence to the parabaloids in ${\Bbb{R}}^3$.

Section 4 turns to reflection in a surface, as expounded more fully in \cite{gak3}. By way of example, we 
compute the reflection of a plane wave off the parabaloids computed in the previous section. In section 5 we 
describe the general procedure for computing the focal set of an arbitrary line congruence and prove Main Theorem 1.

Finally we look in detail at the coffeecup caustics: the focal sets formed by reflection of a point source off
the inside of a cylinder. We do this for multiple reflections off the cylinder and find that the focal sets
consists of two distinct sets: one a translation invariant surface and the other a curve outside of the cylinder (and
hence not physically seen). In the limit as the source tends to infinity, the incoming wave is flat and the resulting 
focal set is independent of the angle of incidence of the wave.  Section 7 contains a discussion of these 
focal sets and their properties.

\section{The Space of Oriented lines}

We start with 3-dimensional Euclidean space ${\Bbb{R}}^3$ and fix standard coordinates
($x^1,x^2,x^3$). In what follows we combine the first two coordinates
to form a single complex coordinate $z=x^1+ix^2$, set $t=x^3$ and refer to coordinates ($z,t$)
on ${\Bbb{R}}^3$.

Let ${\Bbb{L}}$ be the set of oriented lines, or {\it rays}, in Euclidean space
${\Bbb{R}}^3$. Such a line $\gamma$ is uniquely determined by its unit direction vector
$\vec{\mbox{U}}$ and the vector $\vec{\mbox{V}}$ joining the origin to the point on the line
that lies closest to the origin. That is, 

\[
\gamma=\{\;\vec{\mbox{V}}+r\vec{\mbox{U}}\in{\Bbb{R}}^3\;|\;r\in{\Bbb{R}}\;\}
\]
where $r$ is an affine parameter along the line.

By parallel translation, we move $\vec{\mbox{U}}$ to the origin 
and $\vec{\mbox{V}}$ to the head of $\vec{\mbox{U}}$. Thus, we obtain a vector that is tangent to
the unit 2-dimensional sphere in ${\Bbb{R}}^3$. 
The mapping is one-to-one and so it identifies the space of oriented lines 
with the tangent bundle of the 2-sphere ${\mbox{T\:S}}^2$ (see Figure 1).

\[
{\Bbb{L}}=\{\;(\vec{\mbox{U}},\vec{\mbox{V}})\in{\Bbb{R}}^3\times{\Bbb{R}}^3\;
               |\;\quad |\vec{\mbox{U}}|=1\quad\vec{\mbox{U}}\cdot\vec{\mbox{V}}=0\;\}
\]

\vspace{0.1in}
\setlength{\epsfxsize}{4.5in}
\begin{center}
   \mbox{\epsfbox{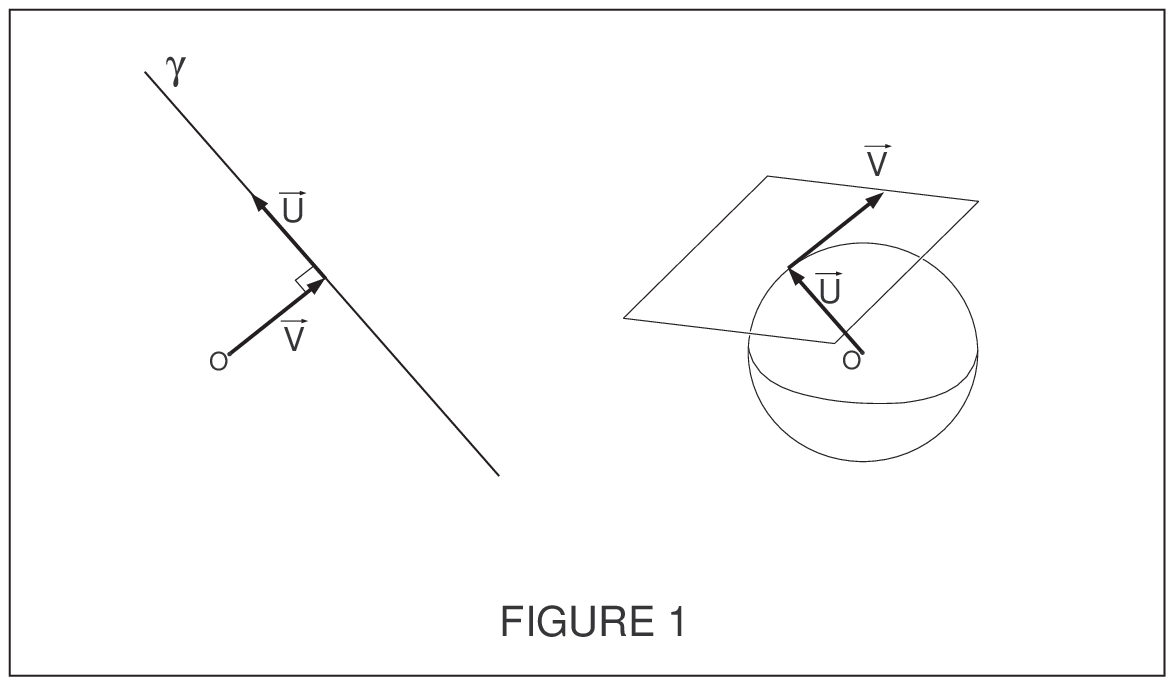}}
\end{center}
\vspace{0.1in}

${\Bbb{L}}$ is a 4-dimensional manifold and the above identification gives a 
natural set of local complex coordinates.
Let $\xi$ be the local complex coordinate on the unit 2-sphere in ${\Bbb{R}}^3$
obtained by stereographic projection from the south pole. 

In terms of the standard spherical polar angles $(\theta,\phi)$, we have
 $\xi=\tan(\frac{\theta}{2})e^{i\phi}$. We convert from coordinates ($\xi,\bar{\xi}$) back to
($\theta,\phi$) using
\begin{equation}\label{e:inv}
\cos\theta={\textstyle{\frac{1-\xi\bar{\xi}}{1+\xi\bar{\xi}}}}
\qquad
\sin\theta={\textstyle{\frac{2\sqrt{\xi\bar{\xi}}}{1+\xi\bar{\xi}}}}
\qquad
\cos\phi={\textstyle{\frac{\xi+\bar{\xi}}{2\sqrt{\xi\bar{\xi}}}}}
\qquad
\sin\phi={\textstyle{\frac{\xi-\bar{\xi}}{2i\sqrt{\xi\bar{\xi}}}}}
\end{equation}

This can be extended to complex
coordinates $(\xi,\eta)$ on ${\Bbb{L}}$ minus the tangent space over the south
pole, as follows. First note that a tangent vector $\vec{\mbox{X}}$ to the 2-sphere can  
always be expressed as a linear combination of the tangent vectors generated by $\theta$ and $\phi$:
\[
\vec{\mbox{X}}=X^\theta\frac{\partial}{\partial\theta}+X^\phi\frac{\partial}{\partial\phi}
\]
In our complex formalism, we have the natural complex tangent vector
\[
\frac{\partial}{\partial\xi}=\cos^2({\textstyle{\frac{\theta}{2}}})
    \left(\frac{\partial}{\partial\theta}
       -\frac{i}{2\cos({\textstyle{\frac{\theta}{2}}})\sin({\textstyle{\frac{\theta}{2}}})}\frac{\partial}{\partial\phi}
    \right) e^{-i\phi}
\]
and any real tangent vector can be written as
\[
\vec{\mbox{X}}=\eta\frac{\partial}{\partial\xi}+\bar{\eta}\frac{\partial}{\partial\bar{\xi}}
\]
for a complex number $\eta$. We identify the real tangent vector $\vec{\mbox{X}}$ on the 2-sphere 
(and hence the ray in ${\Bbb{R}}^3$) with the two complex numbers ($\xi,\eta$). Loosely speaking,
$\xi$ determines the direction of the ray, and $\eta$ determines its perpendicular distance
vector to the origin - complex representations of $\vec{\mbox{U}}$ and $\vec{\mbox{V}}$.

The coordinates ($\xi,\eta$) do not cover all of ${\Bbb{L}}$ - they omit all of the lines 
pointing directly downwards. However, the construction can also be carried out 
using stereographic projection from the north pole, yielding a coordinate system that 
covers all of ${\Bbb{L}}$ except for the lines pointing directly upwards. Between these two 
coordinate patches the whole of the space of oriented lines is covered. In what follows we work in the patch 
that omits the south direction.

While the coordinates depend on a choice of origin, a translation of the origin simply alters the
coordinates by :
\begin{equation}\label{e:trans}
\xi\rightarrow\xi'=\xi  \qquad 
\eta\rightarrow \eta'=\eta+\alpha+b\xi-\bar{\alpha}\xi^2,
\end{equation}
while a rotation about the origin is given by
\begin{equation}\label{e:rot}
\xi\rightarrow\xi'=\frac{a\xi-\beta}{\bar{\beta}\xi+a}  \qquad 
\eta\frac{\partial}{\partial \xi}\rightarrow \eta'\frac{\partial}{\partial \xi'}
               =\frac{\eta}{(\bar{\beta}\xi+a)^2}\frac{\partial}{\partial \xi'},
\end{equation}
for $\alpha,\beta\in{\Bbb{C}}$ and $a,b\in{\Bbb{R}}$ with $a^2+\beta\bar{\beta}=1$. 

\begin{Def}
The map $\Phi:{\Bbb{L}}\times{\Bbb{R}}\rightarrow{\Bbb{R}}^3$ is defined to take 
$((\xi,\eta),r)\in{\Bbb{L}}\times{\Bbb{R}}$ to the point  in ${\Bbb{R}}^3$ on the oriented line ($\xi,\eta$) 
that lies a distance $r$ from the point on the line closest to the origin (see the Figure 2).
\end{Def}

\vspace{0.1in}
\setlength{\epsfxsize}{4.5in}
\begin{center}
   \mbox{\epsfbox{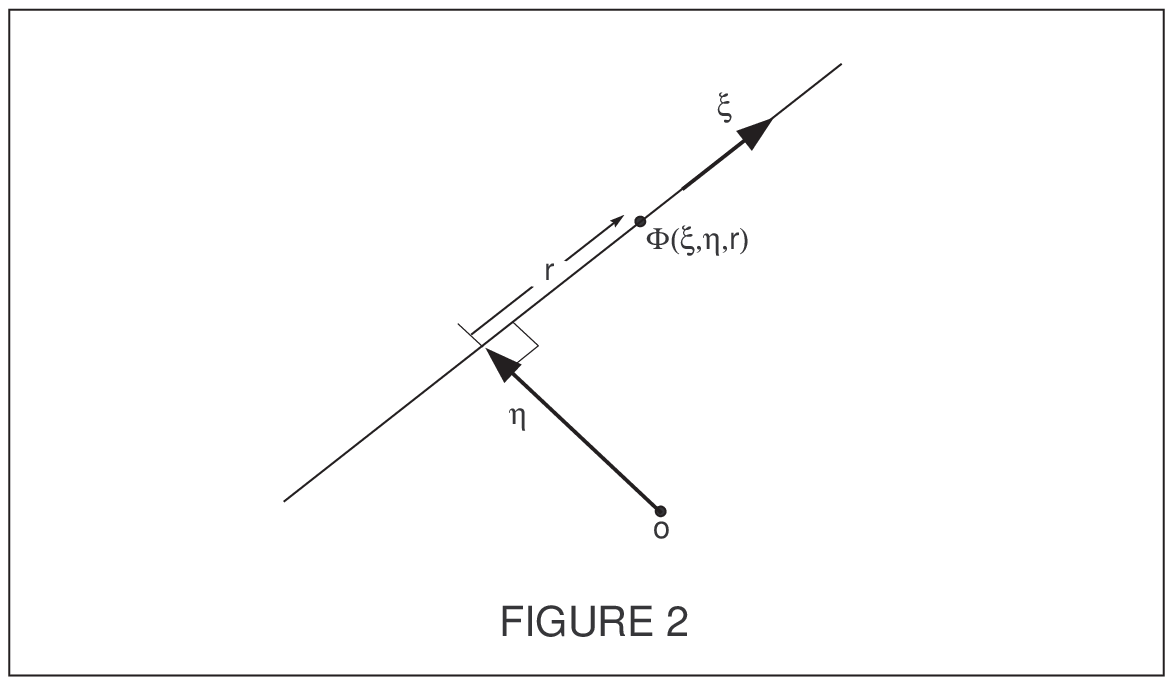}}
\end{center}
\vspace{0.1in}

This map, which is of crucial importance when describing surfaces in ${\Bbb{R}}^3$, 
has the following cooordinate expression:

\begin{Prop}\cite{gak2}
If $\Phi(\xi,\eta,r)=(z(\xi,\eta,r),t(\xi,\eta,r))$, then:
\begin{equation}\label{e:coord}
z=\frac{2(\eta-\overline{\eta}\xi^2)+2\xi(1+\xi\overline{\xi})r}{(1+\xi\overline{\xi})^2}
\qquad\qquad
t=\frac{-2(\eta\overline{\xi}+\overline{\eta}\xi)+(1-\xi^2\overline{\xi}^2)r}{(1+\xi\overline{\xi})^2},
\end{equation}
where $z=x^1+ix^2$, $t=x^3$ and ($x^1$, $x^2$, $x^3$) are Euclidean
coordinates in ${\Bbb{R}}^3$. 
\end{Prop}

The existence of complex coordinates implies that there is a complex structure 
${\Bbb{J}}$ on ${\Bbb{L}}$. In fact, this map ${\Bbb{J}}:{\mbox{T}}_\gamma{\Bbb{L}}\rightarrow {\mbox{T}}_\gamma{\Bbb{L}}$
can be defined as follows: given a line $\gamma$ in ${\Bbb{R}}^3$, the tangent space ${\mbox{T}}_\gamma{\Bbb{L}}$ 
can be identified with
the Jacobi fields along the line that are orthogonal to the direction of $\gamma$. Then ${\Bbb{J}}$ is
given  by rotation of the Jacobi field through 90$^{\mbox{o}}$ about the line $\gamma$ \cite{hitch}. 

The complex structure can be supplemented with a natural symplectic structure $\Omega$. This is obtained by
pulling back the canonical symplectic structure on ${\mbox {T}}^*{\mbox {S}}^2$ with the round metric, 
considered as a mapping from ${\mbox {T\;S}}^2$ to ${\mbox {T}}^*{\mbox {S}}^2$. This symplectic structure
is compatible with the complex structure and so together they define a metric by 
${\Bbb{G}}(\cdot,\cdot)=\Omega({\Bbb{J}}\cdot,\cdot)$.

The metric ${\Bbb{G}}$ has many interesting properties: it is of signature ($++--$), conformally flat and scalar flat
(but not Einstein). In addition, the identity component of the isometry group of ${\Bbb{G}}$ is isomorphic
to the identity component of the Euclidean isometry group \cite{gak4}. A K\"ahler structure with this property
on the space of oriented lines in ${\Bbb{R}}^n$ exists only when $n=3$ or $n=7$ \cite{salvai}. The metric has
the following interpretation: the norm of a vector in ${\mbox{T}}_\gamma{\Bbb{L}}$ with repsect to ${\Bbb{G}}$ 
is the angular momentum of the associated Jacobi field along $\gamma$.

In the next section we investigate the geometric
structures induced by (${\Bbb{J}}$,$\Omega$,${\Bbb{G}}$) on a surface in ${\Bbb{L}}$.

\vspace{0.1in}

\section{Line Congruences}

In its simplest form, geometric optics models the propagation of light through a homogeneous isotropic medium
by a 2-parameter family of rays. 

\begin{Def}
A {\it line congruence} is a 2-parameter family of oriented lines in ${\Bbb{R}}^3$. 
\end{Def}

From our perspective a line congruence is a 
surface $\Sigma$ in ${\Bbb{L}}$. For example, a point source corresponds to the 2-parameter family of 
oriented lines that contain the source point, which thus defines a 2-sphere in ${\Bbb{L}}$.

For computational purposes, we must give explicit local parameterizations of the line congruence. In practice,
it will be given locally by a map ${\Bbb{C}}\rightarrow{\Bbb{L}}:\mu\mapsto(\xi(\mu,\bar{\mu}),\eta(\mu,\bar{\mu}))$. A
convenient choice of parameterization will often depend upon the specifics of the situation, but 
our formalism holds for arbitrary parameterizations.

The dual picture of light propagation is to consider the wavefronts, or surfaces that are orthogonal to a
given set of rays. However, not every line congruence has such orthogonal surfaces - indeed, most don't. To 
explain this we consider the first order properties of $\Sigma$, which can described by two complex functions,
the {\it optical scalars}: $\rho,\sigma:\Sigma\times{\Bbb{R}}\rightarrow{\Bbb{C}}$
The real part $\Theta$ and the imaginary part $\lambda$ of $\rho$ are the {\it divergence} and {\it twist} of the
congruence, while $\sigma$ is the {\it shear}. 

\begin{Def}
A {\it null frame} in ${\Bbb{R}}^3$ is a trio $\{e_0,e_+,e_-\}$ of
complex vector fields in ${\Bbb{C}}\otimes T\:{\Bbb{R}}^3$, where $e_0$ is real, 
$e_+$ is the complex conjugate of $e_-$ and they satisfy the following orthogonality properties:
\[
<e_0, e_0>=1 \qquad <e_0, e_+>=0 \qquad <e_+, e_+>=0 \qquad <e_+,
e_->=1,
\]
where we have extended the Euclidean inner product $<\cdot,\cdot>$ of ${\Bbb{R}}^3$ bilinearly over
${\Bbb{C}}$. Orthonormal frames $\{e_0,e_1,e_2\}$ on $T\:{\Bbb{R}}^3$ and null frames are
related by
\begin{equation}\label{e:ortho}
e_+=\frac{1}{\sqrt{2}}(e_1-ie_2) \qquad
e_-=\frac{1}{\sqrt{2}}(e_1+ie_2).
\end{equation}
\end{Def}

\begin{Def}
A {\it congruence null frame} for $\Sigma\subset{\Bbb{L}}$ is a null frame
$\{e_0,e_+,e_-\}$ if, for each $\gamma\in\Sigma$, we have $e_0$ tangent
to $\gamma$ in ${\Bbb{R}}^3$, and the orientation of $\{e_0,e_1,e_2\}$
is the standard orientation on ${\Bbb{R}}^3$.
\end{Def}

\begin{Def}
Given a line congruence and null frame, the {\it optical scalars} are defined by:
\[
\rho=<\nabla_{e_0}\;e_+,e_->\quad\qquad \sigma=<\nabla_{e_0}\;e_+,e_+>
\]
where $\nabla$ is the Euclidean connection on ${\Bbb{R}}^3$.
\end{Def}

These have the following geometric interpretation.
Consider a specific ray $\gamma$ in the line congruence and a point $p$ along this ray. 
Now consider the unit circle in the plane
orthogonal to the ray at $p$. As we flow this circle along the line congruence this circle will become distorted.  
To first order in the affine parameter along the ray the real part of $\rho$ measures the divergence or 
contraction of the circle, the imaginary part determines the rotation of the circle, while $\sigma$ 
measures the shearing \cite{par} (see Figure 3).

\vspace{0.1in}
\setlength{\epsfxsize}{4.5in}
\begin{center}
   \mbox{\epsfbox{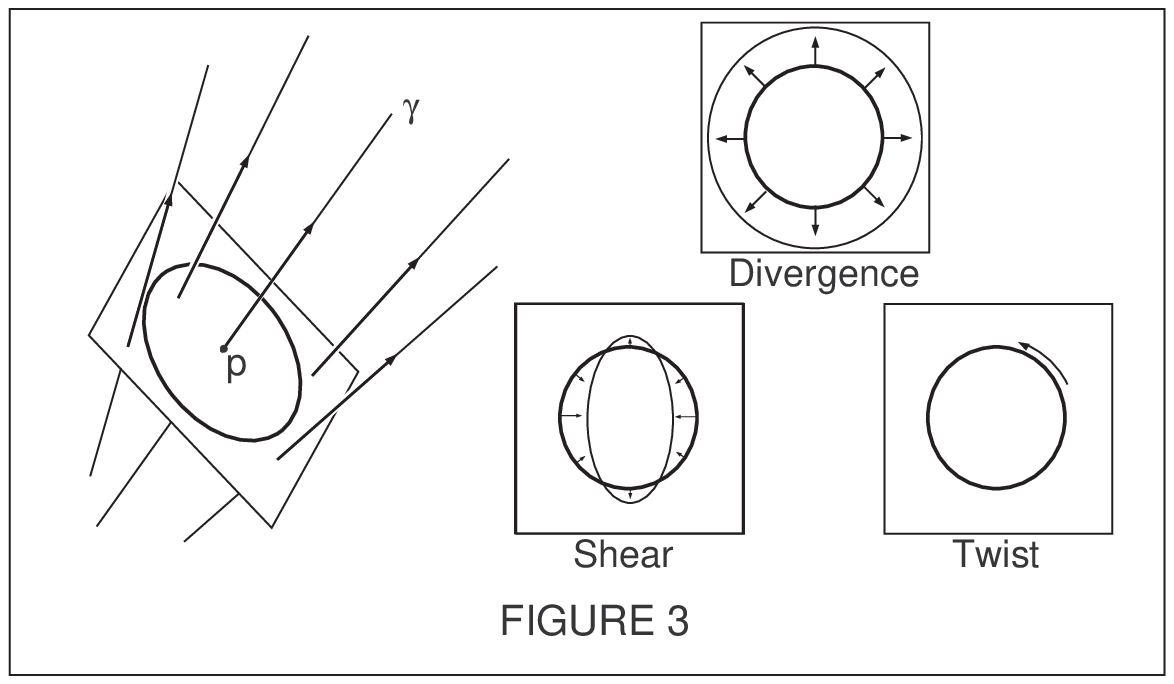}}
\end{center}
\vspace{0.1in}

For computational purposes the optical scalars are given by:

\begin{Prop} \cite{gak2}
For a parameterized line congruence the
optical scalars have the following expressions in terms of first derivatives of the parameterization:  
\begin{equation}\label{e:spinco}
\rho=\theta+\lambda i=\frac{ \partial^+\eta\overline{\partial}\;\overline{\xi} -\partial^-\eta\partial\overline{\xi}}
{\partial^-\eta\overline{\partial^-\eta}-\partial^+\eta\overline{\partial^+\eta}}
\qquad\qquad
\sigma=\frac{\overline{\partial^+\eta}\partial\overline{\xi} -\overline{\partial^-\eta}\;\overline{\partial}\;\overline{\xi}}
{\partial^-\eta\overline{\partial^-\eta}-\partial^+\eta\overline{\partial^+\eta}},
\end{equation}
where
\[
\partial^+\eta\equiv\partial \eta+r\partial\xi-\frac{2\eta\overline{\xi}\partial \xi}{1+\xi\overline{\xi}}
\qquad\qquad
\partial^-\eta\equiv\overline{\partial} \eta+r\overline{\partial}\xi-\frac{2\eta\overline{\xi}\;\overline{\partial} \xi}{1+\xi\overline{\xi}},
\]
and $\partial$ and $\bar{\partial}$ are differentiation with respect
to $\mu$ and $\bar{\mu}$, respectively.
\end{Prop}

The twist has the following important interpretation: 

\begin{Prop}
There exists surfaces orthogonal to the rays of a 
line congruence if and only if the twist of the line congruence is zero. 
\end{Prop}

Moreover, in terms of the K\"ahler structure, the optical scalars have the following significance:

\begin{Thm}\cite{gak4}\label{t:kaehler}

A line congruence $\Sigma\subset{\Bbb{L}}$ is lagrangian (i.e. $\Omega|_\Sigma=0$) iff the twist of $\Sigma$ is zero.

A line congruence $\Sigma\subset{\Bbb{L}}$ is holomorphic (i.e. ${\Bbb{J}}$ preserves the tangent space $T\Sigma$) iff 
the shear of $\Sigma$ is zero.

The metric induced on $\Sigma$ by ${\Bbb{G}}$ is riemannian, degenerate or lorentz iff $|\sigma|^2<\lambda^2$, 
$|\sigma|^2=\lambda^2$ or $|\sigma|^2>\lambda^2$, respectively, where $\lambda$ is the imaginary part of $\rho$.

\end{Thm}

A further geometric quantity is the {\it curvature} of the line congruence, which is defined to be
$\kappa=\rho\bar{\rho}-\sigma\bar{\sigma}$. A line congruence will be said to be {\it flat}
if $\kappa=0$. If a line congruence is non-flat, then the direction of the congruence can be used as a 
parameterization \cite{gak2}. In other words the line congruence is locally given by $\xi\mapsto(\xi,\eta(\xi,\bar{\xi}))$.
The point source line congruence is non-flat, while the set of rays orthogonal to a given line forms a flat line 
congruence.

\section{Constructing Surfaces Using Line Congruences}

We now describe how to construct surfaces in ${\Bbb{R}}^3$ using line
congruences. Given a line congruence $\Sigma\subset{\Bbb{L}}$, a map
$r:\Sigma\rightarrow {\Bbb{R}}$ determines a map
$\Sigma\rightarrow{\Bbb{R}}^3$ by
$(\xi,\eta)\mapsto\Phi((\xi,\eta),r(\xi,\eta))$ for $(\xi,\eta)\in\Sigma$. In other words,
we pick out one point on each line in the congruence (see Figure 4).

With a local parameterization $\mu$ of $\Sigma$, composition with the above map yields a
map ${\Bbb{C}}\rightarrow{\Bbb{R}}^3$ which comes from 
substituting $r=r(\mu,\bar{\mu})$ in equations (\ref{e:coord}).

\vspace{0.1in}
\setlength{\epsfxsize}{4.5in}
\begin{center}
   \mbox{\epsfbox{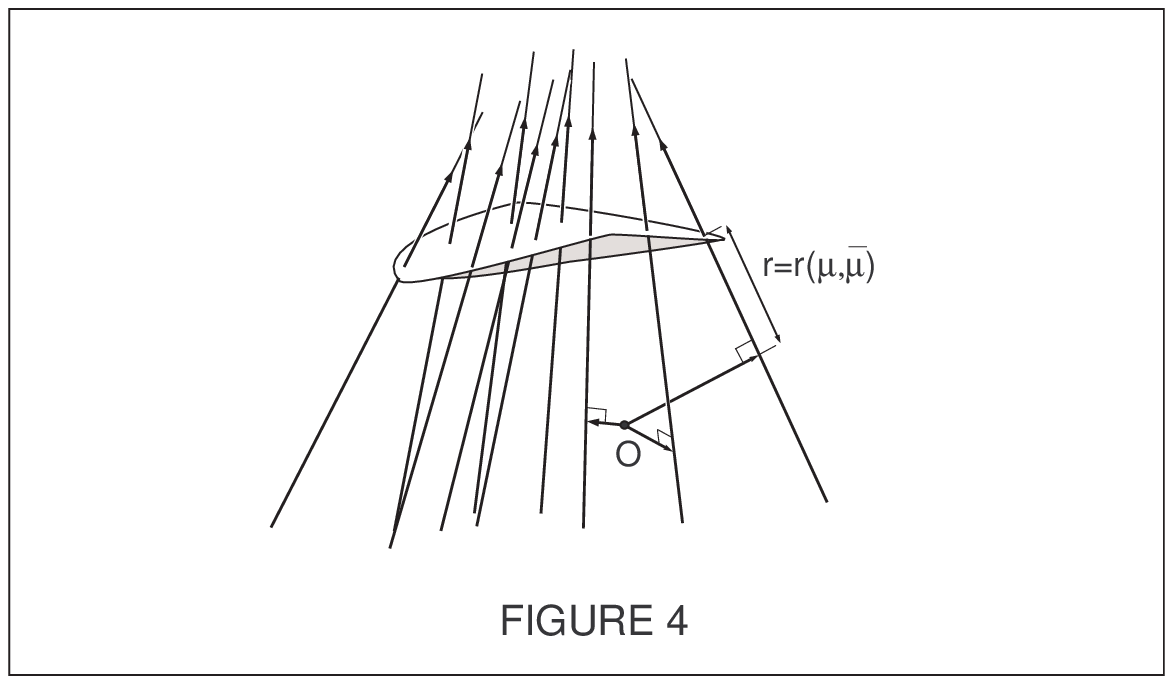}}
\end{center}
\vspace{0.1in}

Of particular interest are the surfaces in ${\Bbb{R}}^3$ orthogonal to
the line congruence - when the line congruence is {\it normal}. As mentioned earlier, these
exist iff the twist of the congruence vanishes. By the first of equation (\ref{e:spinco}), this is 
an integrability condition for a real function:

\begin{Thm}\cite{gak2}

A line congruence $(\xi(\mu,\bar{\mu}),\eta(\mu,\bar{\mu}))$ is orthogonal to a surface in ${\Bbb{R}}^3$ iff there exists a
real function $r(\mu,\bar{\mu})$ satisfying:
\begin{equation}\label{e:intsur}
\bar{\partial} r=\frac{2\eta\bar{\partial}\bar{\xi}+2\bar{\eta}\bar{\partial}\xi}{(1+\xi\bar{\xi})^2}.
\end{equation}
If there exists one solution, there exists a 1-parameter family generated by a real constant of integration.
An explicit parameterization of these surfaces in ${\Bbb{R}}^3$ is given by inserting $r=r(\mu,\bar{\mu})$ in  (\ref{e:coord}). 
\end{Thm}

We now show how to construct the normal line congruence of a non-flat oriented surface $S\subset{\Bbb{R}}^3$. As it is
non-flat, we can parameterize $S$ by the direction $\xi$ of the unit normal. Equivalently, we use the real 
stereographic projection coordinates ($\theta,\phi$).

Suppose $S$ is defined by the pre-image
of zero of a function $G:{\Bbb{R}}^3\rightarrow{\Bbb{R}}$. Then
the unit normal to $S$ is, using standard Euclidean coordinates
($x^1,x^2,x^3$),
\[
\hat{N}=\frac{\mbox{grad} G}{|\mbox{grad}
  G|}=A_1(x^1,x^2,x^3)\frac{\partial}{\partial
  x^1}+A_2(x^1,x^2,x^3)\frac{\partial}{\partial
  x^2}+A_3(x^1,x^2,x^3)\frac{\partial}{\partial x^3}
\]
where
\[
A_i=\left(\sum_{j=1}^n\left(\frac{\partial G}{\partial
      x^j}\right)^2\right)^{-1/2}\;\frac{\partial G}{\partial
      x^i}
\]
The link with the coordinates $\xi=\tan(\frac{\theta}{2})e^{i\phi}$ is
\begin{equation}\label{e:coordtrans1}
A_1(x^1,x^2,x^3)=\cos\phi\sin\theta
\qquad\qquad
A_2(x^1,x^2,x^3)=\sin\phi\sin\theta
\end{equation}
\begin{equation}\label{e:coordtrans2}
A_3(x^1,x^2,x^3)=\cos\theta
\end{equation}
The tangent plane through a
point ($\theta,\phi$) on $S$ is given by
\begin{equation}\label{e:tanplane}
A_1x^1+A_2x^2+A_3x^3=B
\end{equation}
where $A_i$ are given by (\ref{e:coordtrans1}) to
(\ref{e:coordtrans2}) and $B$ is  a function of ($\theta,\phi$)
determined by the surface.

The function $r(\xi,\bar{\xi})$
is the distance of the point on the surface to the point on
the normal line which lies closest to the origin. This
is given by
\[
r=\frac{B}{\sqrt{A_1^2+A_2^2+A_3^2}}=B
\]
Finally, the exact functional relationship between $\xi$ and $\eta$ is given by
\[
\eta(\xi,\bar{\xi})=\frac{1}{2}(1+\xi\bar{\xi})^2\frac{\partial
  r}{\partial \bar{\xi}} 
\]

The task then, reduces to finding $r$ as a function of $\theta$ and
$\phi$, or equivalently, $\xi$ and $\bar{\xi}$. In many simple cases it is possible to solve
(\ref{e:coordtrans1}) to (\ref{e:coordtrans2}) and express
($x^1,x^2,x^3$) as functions of $\theta$ and $\phi$. These can be
directly inserted into (\ref{e:tanplane}) to find $r(\theta,\phi)$.

Once we have $r$ and $\eta$ as functions of $\xi$, equations 
(\ref{e:coord}) give the explicit parameterization of $S$ in terms of $\xi$.

We now work through an example in detail.
\vspace{0.1in}

\noindent {\bf Example}: Elliptic and Hyperbolic Paraboloids

Suppose the surface $S$ is determined by
\[
G=x^3+\frac{(x^1)^2}{a}+\frac{(x^2)^2}{b}=0
\]
for some constants $a$ and $b$. The {\it  elliptic paraboloid} with  $a=1$ $b=1$ and the {\it hyperbolic paraboloid} 
with $a=-b=1$ are graphed below.

The unit normal is
\[
\hat{N}=\left(1+\frac{4(x^1)^2}{a^2}+\frac{4(x^2)^2}{b^2}\right)^{-1/2}\left(\frac{2x^1}{a}\frac{\partial}{\partial x^1}+\frac{2x^2}{b}\frac{\partial}{\partial x^2}+\frac{\partial}{\partial x^3}\right)
\]

We can invert the the relations (\ref{e:coordtrans1}) to
(\ref{e:coordtrans2}) to
\[
x^1=\frac{a}{2}\cos\phi\tan\theta
\qquad
x^2=\frac{b}{2}\cos\phi\tan\theta
\qquad
x^3=-\left(\frac{a}{4}\cos^2\phi+\frac{b}{4}\sin^2\phi\right)\tan^2\theta
\]
Note that the coordinate domain $0\leq\theta<\pi/2$, $0\leq\phi<2\pi$
cover all of the paraboloid.

These give 
\[
B=-\left(a\cos^2\phi+b\sin^2\phi\right)\frac{\sin^2\theta}{4\cos\theta}
\]
Finally converting to holomorphic coordinates to get
\[
r=\frac{a(\xi+\bar{\xi})^2-b(\xi-\bar{\xi})^2}
    {4(1-\xi\bar{\xi})(1+\xi\bar{\xi})}
\qquad
\eta=\frac{a(\xi+\bar{\xi})(1+\xi^3\bar{\xi})
        +b(\xi-\bar{\xi})(1-\xi^3\bar{\xi})}
    {4(1-\xi\bar{\xi})^2}
\]

The equations (\ref{e:coord}) now give the explicit parameterization of the 
paraboloid. Figure 5 shows the resulting parameterization of the parabaloids 
with $a=-b=1$ and $a=b=1$.

\vspace{0.1in}
\setlength{\epsfxsize}{4.5in}
\begin{center}
   \mbox{\epsfbox{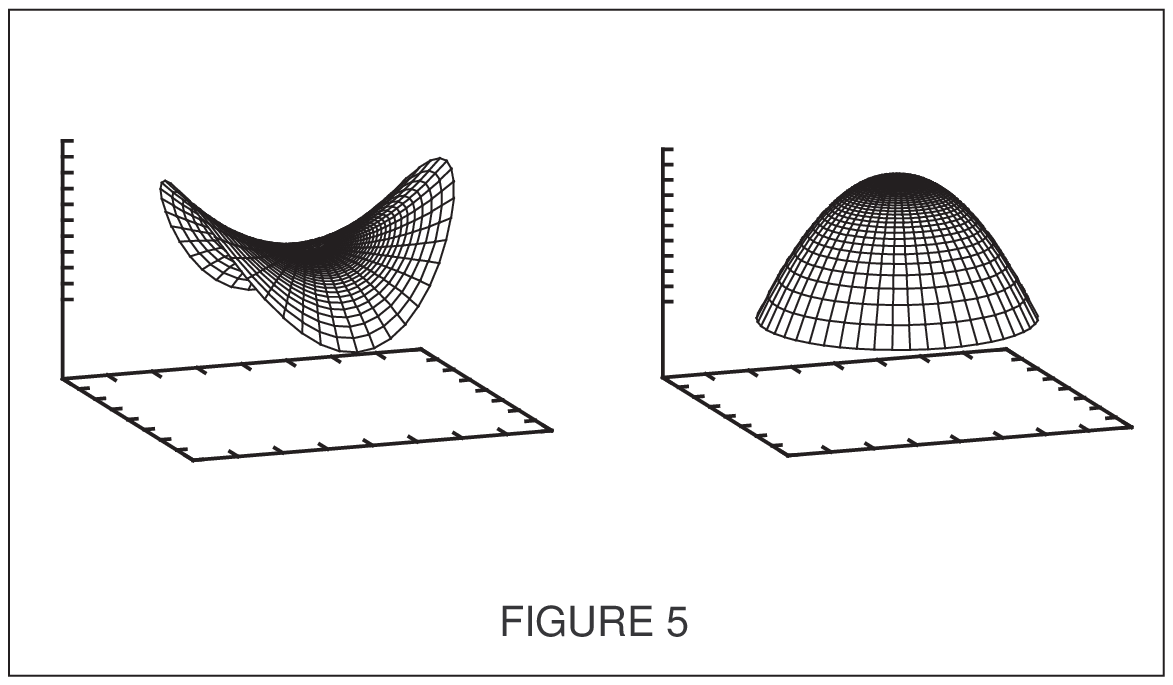}}
\end{center}
\vspace{0.1in}

The lines on these surface are the images of the lines of longitude and
latitude around the north pole under the coordinates.

\section{Reflection}

We turn now to the reflection of an oriented line in a surface in
${\Bbb{R}}^3$. This is equivalent to the action of a certain
group on the space of oriented lines, as described by:

\begin{Thm} \cite{gak3}\label{t:refl}

Consider a parametric line congruence
$\xi=\xi_1(\mu_1,\bar{\mu}_1)$, $\eta=\eta_1(\mu_1,\bar{\mu}_1)$
reflected off an oriented surface with parameterized normal line
congruence $\xi=\xi_0(\mu_0,\bar{\mu}_0)$,
$\eta=\eta_0(\mu_0,\bar{\mu}_0)$ and $r=r_0(\mu_0,\bar{\mu}_0)$
satisfying (\ref{e:intsur}) with $\xi=\xi_0$ and $\eta=\eta_0$. Then
the reflected line congruence ($\xi_2$,$\eta_2$) is 
\begin{equation}\label{e:reflaw1}
\xi_2=\frac{2\xi_0\bar{\xi}_1+1-\xi_0\bar{\xi}_0}
           {(1-\xi_0\bar{\xi}_0)\bar{\xi}_1-2\bar{\xi}_0},
\end{equation}
\begin{equation}\label{e:reflaw2}
\eta_2={\textstyle \frac{(\bar{\xi}_0-\bar{\xi}_1)^2}
         {((1-\xi_0\bar{\xi}_0)\bar{\xi}_1-2\bar{\xi}_0)^2}}\eta_0
       -{\textstyle\frac{(1+\xi_0\bar{\xi}_1)^2}
         {((1-\xi_0\bar{\xi}_0)\bar{\xi}_1-2\bar{\xi}_0)^2}}\bar{\eta}_0
+{\textstyle\frac{(\bar{\xi}_0-\bar{\xi}_1)(1+\xi_0\bar{\xi}_1)(1+\xi_0\bar{\xi}_0)}
         {((1-\xi_0\bar{\xi}_0)\bar{\xi}_1-2\bar{\xi}_0)^2}}r_0,
\end{equation}
where the incoming rays are only reflected if they satisfy the
intersection equation 
\begin{equation}\label{e:int}
\eta_1={\textstyle \frac{(1+\bar{\xi}_0\xi_1)^2}{(1+\xi_0\bar{\xi}_0)^2}}\eta_0
-{\textstyle \frac{(\xi_0-\xi_1)^2}{(1+\xi_0\bar{\xi}_0)^2}}\bar{\eta}_0
+{\textstyle \frac{(\xi_0-\xi_1)(1+\bar{\xi}_0\xi_1)}{1+\xi_0\bar{\xi}_0}}r_0.
\end{equation}
\end{Thm}

\vspace{0.1in}
By virtue of the intersection equation, an alternative way of writing
(\ref{e:reflaw2}) is 
\begin{equation}\label{e:reflaw3}
\eta_2={\textstyle \frac{-(1+\xi_0\bar{\xi}_0)^2}{((1-\xi_0\bar{\xi}_0)\bar{\xi}_1-2\bar{\xi}_0)^2}}\bar{\eta}_1
+{\textstyle \frac{2(\bar{\xi}_0-\bar{\xi}_1)(1+\xi_0\bar{\xi}_1)(1+\xi_0\bar{\xi}_0)}{((1-\xi_0\bar{\xi}_0)\bar{\xi}_1-2\bar{\xi}_0)^2}}r_0.
\end{equation}
The geometric content of this is: reflection of an oriented line
can be decomposed into a sum of a translation and a rotation about the origin.

\vspace{0.1in}

\noindent{\bf Example}: Plane wave reflected off a paraboloid

Consider the paraboloid, as given earlier in section 3. Thus, it is parameterized
by the direction $\xi_0$ of the normal, with
\[
\eta_0=\frac{a(\xi_0+\bar{\xi}_0)(1+\xi_0^3\bar{\xi}_0)
        +b(\xi_0-\bar{\xi}_0)(1-\xi_0^3\bar{\xi}_0)}
    {4(1-\xi_0\bar{\xi}_0)^2}
\]
and
\[
r_0=\frac{a(\xi_0+\bar{\xi}_0)^2-b(\xi_0-\bar{\xi}_0)^2}
    {4(1-\xi_0\bar{\xi}_0)(1+\xi_0\bar{\xi}_0)}
\]

Assume that the incoming plane wave has normal direction along the 
positive $x^1$-axis, that is, $\xi_1=1$. By the reflection equation (\ref{e:reflaw1}), 
the resulting direction is
\[
\xi=\frac{2\xi_0+1-\xi_0\bar{\xi}_0}
           {1-\xi_0\bar{\xi}_0-2\bar{\xi}_0},
\]
and substituting the equation of the paraboloid into (\ref{e:reflaw2}) yields
\begin{align}
\eta=&\frac{a(\xi_0+\bar{\xi}_0)^2(-3-2\xi_0+2\bar{\xi}_0+2\xi_0\bar{\xi}_0
     -2\xi_0\bar{\xi}_0^2+2\xi_0^2\bar{\xi}_0 -3\xi_0^2\bar{\xi}_0^2)}
         {(1-\xi_0\bar{\xi}_0)^2(1-\xi_0\bar{\xi}_0-2\bar{\xi}_0)^2}\\\nonumber
&\qquad+\frac{b(\xi_0-\bar{\xi}_0)(2+3\xi_0-3\bar{\xi}_0+2\xi_0^2+2\bar{\xi}_0^2-2\xi_0\bar{\xi}_0
     +2\xi_0\bar{\xi}_0^2-2\xi_0^2\bar{\xi}_0)}
         {(1-\xi_0\bar{\xi}_0)^2(1-\xi_0\bar{\xi}_0-2\bar{\xi}_0)^2} \\\nonumber
&\qquad\qquad+\frac{b(\xi_0-\bar{\xi}_0)(2\xi_0^2\bar{\xi}_0^2-2\xi_0\bar{\xi}_0^3
     -2\xi_0^3\bar{\xi}_0-3\xi_0^2\bar{\xi}_0^3+3\xi_0^3\bar{\xi}_0^2-2\xi_0^3\bar{\xi}_0^3)}
         {(1-\xi_0\bar{\xi}_0)^2(1-\xi_0\bar{\xi}_0-2\bar{\xi}_0)^2}\nonumber
\end{align}

Since the incoming ray direction is fixed, we are parameterizing the reflected line congruence by the direction 
$\xi_0$ of the normal to the surface at the point of reflection. A direct integration of 
equation (\ref{e:intsur}) gives the function $r$ as
\[
r=-\frac{2(\xi_0+\bar{\xi}_0)[a(\xi_0+\bar{\xi}_0)^2-b(\xi_0-\bar{\xi}_0)^2]}
      {(1+\xi_0\bar{\xi}_0)^2(1-\xi_0\bar{\xi}_0)}+C
\]

Finally, the wavefronts form a one-parameter family
of parameterized surfaces, which can be obtained by substituting for $\xi$, $\eta$ and $r$ in 
equation (\ref{e:coord}). The result in spherical polar coordinates, after some simplifications,
boils down to:
\[
x^1=8a(1-\cos\theta)\tan\theta\cos^3\phi+C (1-2\sin^2\theta\cos^2\phi)
\]
\[
x^2=2(2a\sin^2\theta\cos^2\phi+b)\tan\theta\sin\phi-2C \sin^2\theta\sin\phi\cos\phi
\]
\[
x^3=a(4\cos^2\theta-1)\tan^2\theta\cos^2\phi-b\tan^2\theta\sin^2\phi-2C \sin\theta\cos\theta\cos\phi
\]
where $0\le\theta<\pi/2$ and $0\le\phi<2\pi$.

\vspace{0.1in}
\setlength{\epsfxsize}{4.5in}
\begin{center}
   \mbox{\epsfbox{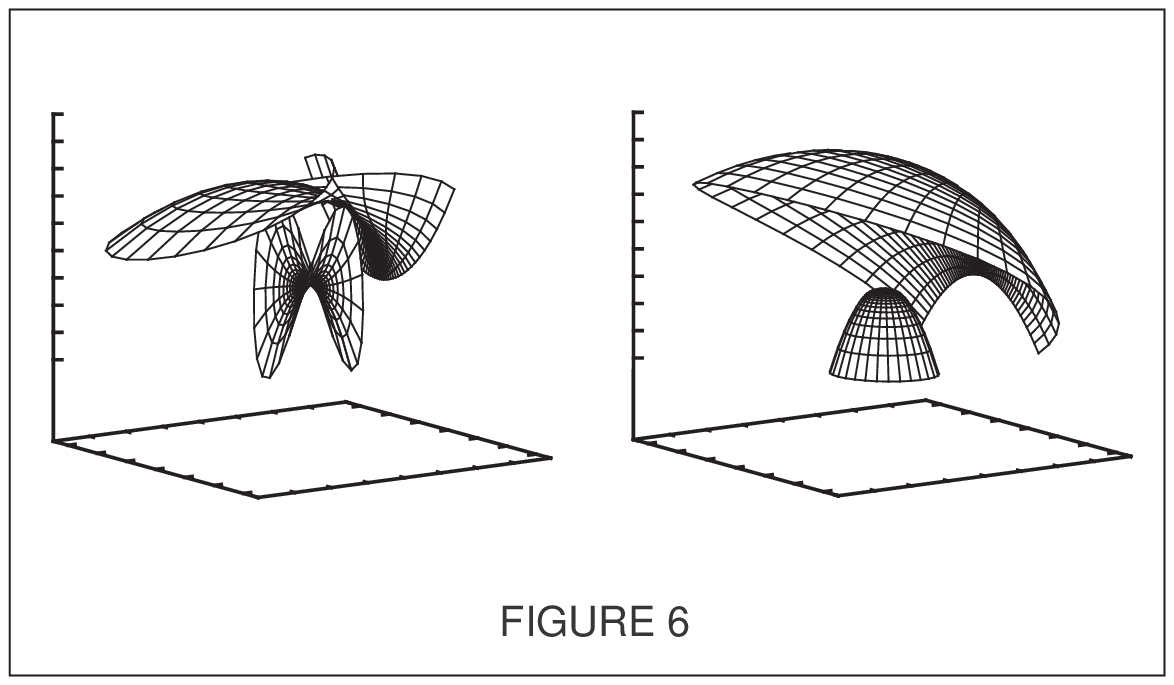}}
\end{center}
\vspace{0.1in}

Figure 6 shows examples of the reflected wavefronts for an 
elliptic and a hyperbolic paraboloid.

\section{The Focal Set of a Line Congruence}

Let $\sigma$ and $\rho$ be the optical scalars of a line congruence
$\Sigma$ as described in section 2. 

\begin{Def}
A point $p$ on a line $\gamma$ in a line congruence is a
{\it focal point} if $\rho$ and $\sigma$ blow-up at $p$. The set of
focal points of a line congruence $\Sigma$ generically form surfaces
in ${\Bbb{R}}^3$, which will be referred to as the {\it focal
surfaces} of $\Sigma$. 
\end{Def}

\begin{Thm}\label{t:focs}

The focal set of a parametric line congruence $\Sigma$ is
\[
\{\;\Phi(\gamma,r)\;|\;\gamma\in\Sigma\;\;\mbox{and}\;\;
1-(\rho_0+\bar{\rho}_0)r+(\rho_0\overline{\rho}_0-\sigma_0\overline{\sigma}_0)r^2=0\}, 
\]
where the coefficients of the quadratic equation are given locally by
(\ref{e:spinco}) evaluated at $r=0$. 
\end{Thm}
\begin{pf}
In terms of the affine
parameter $r$ along a given line, the Sachs equations, which $\sigma$ and
$\rho$ must satisfy, are \cite{par}:
\[
\frac{\partial \rho}{\partial r}=\rho^2+\sigma\overline{\sigma}
\qquad 
\frac{\partial \sigma}{\partial r}=(\rho+\bar{\rho})\sigma.
\]
These are equivalent to the vanishing of certain components of the
Ricci tensor of the Euclidean metric. They have solution:
\[
\rho=\frac{\rho_0-(\rho_0\overline{\rho}_0-\sigma_0\overline{\sigma}_0)r}
  {1-(\rho_0+\bar{\rho}_0)r+(\rho_0\overline{\rho}_0
                              -\sigma_0\overline{\sigma}_0)r^2}
\qquad
\sigma=\frac{\sigma_0}
  {1-(\rho_0+\bar{\rho}_0)r+(\rho_0\overline{\rho}_0-\sigma_0\overline{\sigma}_0)r^2},
\]
where  $\sigma_0$ and $\rho_0$ are the values of the optical
scalars at $r=0$. The theorem follows.
\end{pf}

This has the following corollary:

\begin{Cor}\label{c:foc}
Let $\Sigma$ be a line congruence, $\rho=\Theta+\lambda i$,
$\sigma$ the associated optical scalars and $\rho_0$, $\Theta_0$,
$\lambda_0$, $\sigma_0$ their values at $r=0$. 

If $\Sigma$ is flat with non-zero divergence, then there
exists a unique focal surface $S$ given by
$r=(2\Theta_0)^{-1}$. If it is flat with zero divergence,
then the focal set is empty.

If $\Sigma$ is non-flat, then there exists a unique focal point on each
line iff  $|\sigma_0|^2=\lambda_0^2$, there exist two focal
points on each line iff $|\sigma_0|^2<\lambda_0^2$ and there are no
focal points on each line iff $|\sigma_0|^2>\lambda_0^2$. The
focal set is given by
\[
r=\frac{\Theta_0\pm(|\sigma_0|^2 -\lambda_0^2)^{\frac{1}{2}}}
            {\rho_0\bar{\rho}_0-\sigma_0\bar{\sigma}_0}.
\]

\end{Cor}

\begin{pf}
The focal set of a parameterized line congruence
are given by $r=r(\mu,\bar{\mu})$ satisfying the quadratic equation
in Theorem \ref{t:focs}. If $\kappa=0$, then there is none or one
solution depending on whether $\Theta_0=0$ or not.

If $\kappa\neq 0$ then there are two, one or no solutions iff
$|\sigma_0|^2-\lambda_0^2$ is greater than, equal to or less than
zero (respectively). 

The solution of the quadratic equation in each case is as stated.
\end{pf}

There is also the equivalent definition for focal surfaces:

\begin{Prop}\label{p:2defs}
A continuously differentiable surface $S$ in ${\Bbb{R}}^3$ is a focal surface of a line
congruence $\Sigma$ iff every line in $\Sigma$ is tangent to $S$ at
some point. 
\end{Prop}
\begin{pf}
Let $\Sigma$ be locally parameterized by $\mu\mapsto
(\xi(\mu,\bar{\mu}),\eta(\mu,\bar{\mu}))$, where $(\xi,\eta)$ are the
canonical coordinates above. Then a surface in ${\Bbb{R}}^3$ given
by $r:\Sigma\rightarrow{\Bbb{R}}$ is tangent to the line congruence iff
\[
\mbox{Det}\left[\begin{array}{ccc}
\frac{2\xi}{1+\xi\bar{\xi}}
&\frac{2\bar{\xi}}{1+\xi\bar{\xi}}
&\frac{1-\xi\bar{\xi}}{1+\xi\bar{\xi}}\\
\partial z & \partial\bar{z} & \partial t
\\
\bar{\partial} z & \bar{\partial}\bar{z} & \bar{\partial} t
\end{array}\right]=0,
\]
where the partial derivatives are in $\mu$ and $\bar{\mu}$. This 
determinant equation, which is
\[
2\xi(\partial\bar{z}\bar{\partial} t-\bar{\partial}\bar{z}\partial t)
+2\bar{\xi}(\partial z\bar{\partial} t-\bar{\partial} z\partial t)
+(1-\xi\bar{\xi})(\partial z\bar{\partial}\bar{z}-\bar{\partial} z\partial\bar{z})=0,
\]
is a quadratic equation for $r=r(\mu,\bar{\mu})$ with coefficients
given by the first derivatives of $z(\mu,\bar{\mu})$,
$\bar{z}(\mu,\bar{\mu})$ and $t(\mu,\bar{\mu})$. 

Carrying out the differentiation we find, for example, that
\begin{align}
\partial z&= {\textstyle \left(\frac{2r}{(1+\xi\bar{\xi})^2}
              -\frac{4(\bar{\xi}\eta+\xi\bar{\eta})}
                 {(1+\xi\bar{\xi})^3}\right)}\partial\xi   
+{\textstyle \left(-\frac{2(1+2\xi\bar{\xi})\xi^2r}{(1+\xi\bar{\xi})^2}
              -\frac{4(\eta-\xi^2\bar{\eta})}
                 {(1+\xi\bar{\xi})^3}\right)}\partial\bar{\xi}\nonumber\\
&\qquad\qquad\qquad
+{\textstyle \frac{2}{(1+\xi\bar{\xi})^2}}\partial\eta   
   -{\textstyle \frac{2\xi^2}{(1+\xi\bar{\xi})^2}}\partial\bar{\eta}
   +{\textstyle \frac{2\xi}{1+\xi\bar{\xi}}}\partial r .\nonumber
\end{align}
Similar computations finally yield the quadratic that appears in
Theorem \ref{t:focs}.
\end{pf}

We now explore the geometric properties of the focal set. First we prove:

\vspace{0.1in}

\noindent {\bf Main Theorem 1.}

{\it
Let $\Sigma$ be an immersed surface in ${\Bbb{L}}$. If $\Sigma$ is flat, there is exactly one focal point on each line of
the congruence. If $\Sigma$ is not flat then there is none, one or two focal points on each line iff the
metric induced on $\Sigma$ by ${\Bbb{G}}$ is riemannian, degenerate or lorentz (respectively).
}

\begin{pf}
The number of focal points on a given line is determined by the sign of the discriminant of the quadratic
equation in Theorem \ref{t:focs}: $|\sigma_0|^2-\lambda_0^2$. By Theorem \ref{t:kaehler} ({\it cf.} 
Theorem 2 of \cite{gak4}), this is precisely what determines the sign of the metric induced on $\Sigma$ by ${\Bbb{G}}$: 
the metric is riemannian, degenerate or lorentz iff $|\sigma_0|^2-\lambda_0^2$ is less than, equal to or 
greater than zero. 

The result follows.
\end{pf}

\vspace{0.1in}

Consider now the case where there are two focal points on each line of $\Sigma$. Thus $|\sigma_0|^2-\lambda_0^2>0$,
and further suppose that these focal points form two continuously differentiable surfaces $S_1$ and $S_2$ in ${\Bbb{R}}^3$.
Let $L$ be the distance between the focal points and $\varphi$ the angle between the normals to $S_1$ and $S_2$ at
corresponding points.

\begin{Thm}

The distance $L$ and angle $\varphi$ defined above are given by
\[
L=2\frac{\left(|\sigma_0|^2-\lambda_0^2\right)^{\scriptstyle{\frac{1}{2}}}}{\rho_0\bar{\rho}_0-\sigma_0\bar{\sigma}_0}
\qquad\qquad
\cos^2\varphi=\frac{\lambda_0^2}{|\sigma_0|^2}
\]
\end{Thm}

\begin{pf}
The first of these follows trivially from the fact that the two focal surfaces are given by 
\[
r_1=\frac{\Theta_0+(|\sigma_0|^2 -\lambda_0^2)^{\frac{1}{2}}}
            {\rho_0\bar{\rho}_0-\sigma_0\bar{\sigma}_0}
\qquad\qquad
r_2=\frac{\Theta_0-(|\sigma_0|^2 -\lambda_0^2)^{\frac{1}{2}}}
            {\rho_0\bar{\rho}_0-\sigma_0\bar{\sigma}_0}
\]
The line congruence $\Sigma$, by assumption, is not flat, and so we parameterize it by its direction $\xi$.
To compute the angle $\varphi$ we note that parametric equations for $S_1$ and $S_2$ 
\[
z=z_1(\xi,\bar{\xi})\quad t=t_1(\xi,\bar{\xi}) \qquad\qquad z=z_2(\xi,\bar{\xi})\quad t=t_2(\xi,\bar{\xi})
\]
are obtained by inserting $r=r_1$ and $r=r_2$ in equations (\ref{e:coord}). Let $\nu_1,\nu_2\in{\Bbb{P}}^1$
be the directions of the normals to $S_1$ and $S_2$, respectively. Thus, for $i=1,2$,
\[
\frac{\nu_i}{1-\nu_i\bar{\nu}_i}=\frac{\partial z_i\bar{\partial} t_i-\partial t_i\bar{\partial} z_i}
   {\partial \bar{z}_i\bar{\partial} z_i-\partial z_i\bar{\partial} \bar{z}_i}
\]
If we introduce, for $i=1,2$,
\[
\alpha_i=\partial z_i\bar{\partial} t_i-\partial t_i\bar{\partial} z_i
\qquad\qquad
b_i=\partial \bar{z}_i\bar{\partial} z_i-\partial z_i\bar{\partial} \bar{z}_i
\]
a straightforward computation shows that
\[
\cos\varphi=\pm\frac{b_1b_2+2(\alpha_1\bar{\alpha}_2+\alpha_2\bar{\alpha}_1)}
    {\left[(b_1^2+4\alpha_1\bar{\alpha}_1)(b_2^2+4\alpha_2\bar{\alpha}_2)\right]^{\scriptstyle{\frac{1}{2}}}}
\]
A lengthy computation involving the explicit expressions for $\alpha_i$ and $b_i$ obtained by differentiation of  
(\ref{e:coord}), yields
\begin{align}
b_1b_2+2(\alpha_1\bar{\alpha}_2+\alpha_2\bar{\alpha}_1)&= \frac{4\lambda_0i}{\kappa_0^2(1+\xi\bar{\xi})^2}\left[
  \bar{\sigma}_0(\partial L)^2-\sigma_0(\bar{\partial} L)^2+2\lambda_0i\partial L\bar{\partial} L\right.\nonumber\\
& \qquad\left.+2L(\bar{\beta}\partial L-\beta\bar{\partial} L)-4\beta^2\bar{\sigma}_0+4\bar{\beta}^2\sigma_0
  -8\lambda_0i\beta\bar{\beta}\right]\nonumber
\end{align}
and
\begin{align}
(b_1^2+4\alpha_1\bar{\alpha}_1)(b_2^2+4\alpha_2\bar{\alpha}_2)&= -\frac{16\sigma_0\bar{\sigma}_0}
    {\kappa_0^4(1+\xi\bar{\xi})^4}\left[
  \bar{\sigma}_0(\partial L)^2-\sigma_0(\bar{\partial} L)^2+2\lambda_0i\partial L\bar{\partial} L\right.\nonumber\\
& \qquad\left.+2L(\bar{\beta}\partial L-\beta\bar{\partial} L)-4\beta^2\bar{\sigma}_0+4\bar{\beta}^2\sigma_0
  -8\lambda_0i\beta\bar{\beta}\right]^2\nonumber
\end{align}
where we have introduced 
\[
\beta=(1+\xi\bar{\xi})^2\bar{\partial}\left(\frac{\sigma_0}{\kappa_0(1+\xi\bar{\xi})^2}\right)
     +i\partial\left(\frac{\lambda_0}{\kappa_0}\right).
\]
The expression for the angle $\varphi$ follows.
\end{pf}

\section{Reflection off a Cylinder}

Consider a cylinder of radius $a$, with axis lying along
the $x^3-$axis in ${\Bbb{R}}^3$. Then we have:

\begin{Prop}
The inward pointing normal to such a cylinder is
given parametrically by:
\[
\xi=e^{iv}\qquad\qquad \eta=-u\;e^{iv} 
\]
for $(u,v)\in{\Bbb{R}}\times\mbox{S}^1$. The distance of a point $p$ on the surface from the
point on the normal through $p$ that lies closest to the origin is $r=-a$.
\end{Prop}
\begin{pf}
This can be checked by noting that, with the aid of (\ref{e:coord}), the mapping 
$(u,v)\mapsto\Phi(\xi(u,v),\eta(u,v),r(u,v))$, with $\xi$, $\eta$ and $r$ as stated, 
yields a parameterization of the cylinder:
$(u,v)\mapsto(-a\cos v,-a\sin v,u)$. Moreover, the oriented normal at the point $(u,v)$ on the cylinder is given by 
the expression in the proposition. 
\end{pf}

The coffeecup caustic is obtained by finding the focal set of a plane
wave reflected off the inside of this cylinder. This turns out to be:

\begin{Prop}\label{p:cylpl}
Consider the reflection off the inside of a cylinder of
radius $a$ of a line congruence consisting of parallel rays traveling along the $x^1$-axis making an
angle $\beta$ with the $x^3$-axis. The focal set of the reflected line
congruence is a surface given parametrically by 
\[
x^1=a\cos v\left(\cos^2v-\frac{3}{2}\right)
\qquad\qquad
x^2=a\sin v\left(\cos^2v-1\right),
\]
\[
x^3=-u-\frac{a}{2}\cos v\cot\beta,
\]
for $u\in{\Bbb{R}}$ and $\pi/2\le v\le3\pi/2$.
\end{Prop}

\begin{pf}
The normal congruence to the plane wave is ($\xi_1$,$\eta_1$), where $\xi_1$ and $\eta_1\in{\Bbb{C}}$ is free. 
Reflecting this off the cylinder given above, we have by Theorem \ref{t:refl}:
\[
\xi_2=-\bar{\xi}_1e^{2iv}
\qquad\qquad
\eta_2=-{\textstyle{\frac{1}{2}}}\left(ae^{-iv}-2u\bar{\xi}_1-ae^{iv}\bar{\xi}_1^2 \right)e^{2iv}
\]
or if we parameterize by the point of reflection $\mu=\eta_0=-ue^{iv}$
\[
\xi_2=-\frac{\mu\bar{\xi}_1}{\bar{\mu}}
\qquad\qquad
\eta_2=-{\textstyle{\frac{1}{2}}}\left(ae^{-iv}-2\mu\bar{\xi}_1
-a{\textstyle{\frac{\mu}{\bar{\mu}}}}\bar{\xi}_1^2 \right)\left({\textstyle{\frac{\mu}{\bar{\mu}}}}\right)
                ^{\scriptstyle{\frac{1}{2}}}
\]
We now compute the optical scalars for this line congruence via equations (\ref{e:spinco}) and find that the
line congruence is flat (i.e. $\rho\bar{\rho}-\sigma\bar{\sigma}=0$). Thus by Corollary \ref{c:foc} there is exactly
one focal point on each line, given by, after some computation:
\[
r=-u\cos\beta-\frac{a\cos v(2\cos^2\beta-1)}{2\sin\beta}
\] 
Inserting this, along with the expressions for $\xi_2$ and $\eta_2$ in (\ref{e:coord}) yields the stated result.
The domain of $v$ must be restricted to half a circle as the incoming rays reflect on the inside of only one half of the
cylinder.
\end{pf}

We now consider the focal set formed by reflection of a point
source off the inside of the cylinder. To this end, the following theorem describes reflection in a cylinder
as a mapping ($\xi_1$,$\eta_1$)$\mapsto$($\xi_2$,$\eta_2$):

\begin{Thm}\label{t:cylref}

A ray ($\xi_1$, $\eta_1$) intersects a
cylinder of radius $a$ lying along the $x^3$-axis iff:
\begin{equation}\label{e:cylint1}
\left|\frac{\xi_1\bar{\eta}_1-\bar{\xi}_1\eta_1}
   {\xi_1(1+\xi_1\bar{\xi}_1)}\right|\le a.
\end{equation}
For such a ray, the reflected ray is
\[
\xi_2=-\bar{\xi}_1\left( \frac{\xi_1\bar{\eta}_1-\bar{\xi}_1\eta_1
    \pm(a^2\xi_1\bar{\xi}_1(1+\xi_1\bar{\xi}_1)^2
           +(\xi_1\bar{\eta}_1-\bar{\xi}_1\eta_1)^2)^{\frac{1}{2}}}
      {a\bar{\xi}_1(1+\xi_1\bar{\xi}_1)}\right)^2,
\]
\begin{align}
\eta_2=&-\frac{1}{\xi_1}\left(\bar{\xi}_1\eta_1
          \pm\frac{1-\xi_1\bar{\xi}_1}{1+\xi_1\bar{\xi}_1}
          (a^2\xi_1\bar{\xi}_1(1+\xi_1\bar{\xi}_1)^2
           +(\xi_1\bar{\eta}_1-\bar{\xi}_1\eta_1)^2)^{\frac{1}{2}}\right)\nonumber\\
&\qquad\qquad\qquad.
\left(\frac{\xi_1\bar{\eta}_1-\bar{\xi}_1\eta_1
    \pm(a^2\xi_1\bar{\xi}_1(1+\xi_1\bar{\xi}_1)^2
           +(\xi_1\bar{\eta}_1-\bar{\xi}_1\eta_1)^2)^{\frac{1}{2}}}
      {a\bar{\xi}_1(1+\xi_1\bar{\xi}_1)}\right)^2,
\end{align}
where exterior and interior reflection are given by the plus and minus
signs, respectively.
\end{Thm}
\begin{pf}
Consider an incoming ray ($\xi_1$, $\eta_1$). The
reflection equations (\ref{e:reflaw1}) and (\ref{e:reflaw2}) tell us again
that the reflected ray is 
\begin{equation}\label{e:refincyl}
\xi_2=-\bar{\xi}_1e^{2iv}\qquad\qquad
  \eta_2=-\frac{1}{2}\left(ae^{-iv}-2u\bar{\xi}_1
    -ae^{iv}\bar{\xi}_1^2\right)e^{2iv}.
\end{equation}
This incoming ray intersects the cylinder iff ({\it cf} (\ref{e:int})): 
\begin{equation}\label{e:cylint2}
\eta_1=-\frac{1}{2}\left(ae^{iv}+2u\xi_1
    -ae^{-iv}\xi_1^2\right).
\end{equation}
We eliminate $u$ from this equation by combining it with its conjugate
and solving the resulting equation for $v$. The solution, which
exists iff (\ref{e:cylint1}) holds, is
\[
e^{iv}= \frac{\xi_1\bar{\eta}_1-\bar{\xi}_1\eta_1
    \pm(a^2\xi_1\bar{\xi}_1(1+\xi_1\bar{\xi}_1)^2
           +(\xi_1\bar{\eta}_1-\bar{\xi}_1\eta_1)^2)^{\frac{1}{2}}}
      {a\bar{\xi}_1(1+\xi_1\bar{\xi}_1)}.
\]
Substituting this back into the intersection equation we get that
\[
u=-\frac{1}{2\xi_1\bar{\xi}_1}\left({\textstyle \xi_1\bar{\eta}_1+\bar{\xi}_1\eta_1
          \pm\frac{1-\xi_1\bar{\xi}_1}{1+\xi_1\bar{\xi}_1}
          (a^2\xi_1\bar{\xi}_1(1+\xi_1\bar{\xi}_1)^2
           +(\xi_1\bar{\eta}_1-\bar{\xi}_1\eta_1)^2)^{\frac{1}{2}}}\right).
\]
Finally, putting these last two equations into the reflected ray equation
(\ref{e:refincyl}) yields the stated result.

\end{pf}

For multiple reflection we have the following:

\begin{Thm}\label{t:kthreflcyl}

The k$^{th}$ reflection of a ray ($\xi_1$, $\eta_1$) off the inside of
the  cylinder is:
\[
\xi_{k+1}=(-1)^k\xi_1\left[\frac{\Psi_1i-(|\xi_1|^2-\Psi_1^2)^{\frac{1}{2}}}{|\xi_1|}\right]^{2k}
\]
\[
\eta_{k+1}=\frac{(-1)^k}{\bar{\xi}_1}\left[\bar{\xi}_1\eta_1-ka(1-|\xi_1|^2)(|\xi_1|^2-\Psi_1^2)^{\frac{1}{2}}\right]
  \left[\frac{\Psi_1i-(|\xi_1|^2-\Psi_1^2)^{\frac{1}{2}}}{|\xi_1|}\right]^{2k}
\]
where 
\[
\Psi_1=\frac{\xi_1\bar{\eta}_1-\bar{\xi}_1\eta_1}{ai(1+|\xi_1|^2)}
\]
\end{Thm}
\begin{pf}
This follows from iterations of the above theorem once we realize that $|\xi_l|$ and 
\[
\Psi_l\equiv \frac{\xi_l\bar{\eta}_l-\bar{\xi}_l\eta_l}{a(1+|\xi_l|^2)}
\]
are preserved by reflection in a cylinder.
\end{pf}

For a point source at a finite distance, the following theorem
describes the focal set of the $k^{\mbox{th}}$ reflection:

\vspace{0.1in}

\noindent {\bf Main Theorem 2.}

{\it
Consider the k$^{th}$ reflection off the inside of a cylinder lying along x$^3$-axis with
radius $a$ of a point source at ($-l,0,0$).
The focal set of the reflected line congruence consists of 
a surface:
\[
z=(-1)^{k+1}l\frac{[l\sin v\;i+(a^2-l^2\sin^2v)^{\frac{1}{2}}]^{2k}
     [2kl\cos v\sin v\;e^{-iv}+(a^2-l^2\sin^2v)^{\frac{1}{2}}]}
     {a^{2k}[2kl\cos v+(a^2-l^2\sin^2v)^{\frac{1}{2}}]}
\]
\[
x^3=\frac{k(1-u^2)[a^2-l^2-2l^2\sin^2v+2kl\cos v(a^2-l^2\sin^2v)^{\frac{1}{2}}]}{u[2kl\cos v+(a^2-l^2\sin^2v)^{\frac{1}{2}}]}
\]
and a curve in the $x^1x^2$-plane: 
\[
z=(-1)^{k+1}ka^{-2k}[l\sin v\;i+(a^2-l^2\sin^2v)^{\frac{1}{2}}]^{2k}
  [l+2ke^{iv}(a^2-l^2\sin^2v)^{\frac{1}{2}}]
\]
where $z=x^1+ix^2$, $u\in{\Bbb{R}}$ and $v$ is in the domain
\[
0\le v\le\pi \qquad\qquad\mbox{for} \qquad l\le a 
\]
and
\[-\sin^{-1}(a/l)\le v\le\sin^{-1}(a/l)
\qquad\qquad\mbox{for}\qquad l> a.
\]
}

\begin{pf}
Consider a point source lying at ($-l$, $0$, $0$) in
${\Bbb{R}}^3$. This line congruence can be parameterized by its
direction $\xi_1$ and $\eta_1=-l(1-\xi_1^2)/2$.
The line congruence obtained from $k$ reflections of this point source off the
inside of a cylinder of radius $a$ is given by Theorem \ref{t:kthreflcyl}.

We then compute the optical scalars of this line congruence parameterized by $\xi_1$ using (\ref{e:spinco}). 
We find, for example, that
\[
\frac{|\sigma_0|}{\kappa_0}=\frac{k(1+u^2)\left(kl\cos v+(a^2-l^2\sin^2v)^{\frac{1}{2}}\right)^2}
           {2u\left(2kl\cos v+(a^2-l^2\sin^2v)^{\frac{1}{2}}\right)},
\]
where $\xi_1=ue^{iv}$. The similar expression for $\rho_0$ shows that
the reflected congruence is not flat. Thus each line contains exactly
two focal points which can be obtained by inserting the solutions of
the quadratic equation of Theorem \ref{t:focs} into
(\ref{e:coord}). The results are as stated above. 
\end{pf}

\section{Discussion}

The focal surface obtained from a plane wave reflected off the inside of a cylinder (Proposition
\ref{p:cylpl}) is symmetric along the $x^3$-axis and intersects
any plane parallel to the $x^1x^2-$plane in a curve. This curve, called a nephroid, is often observed on 
the top of a cup of coffee in the presence of a strong, low and distant light source - hence the
sobriquet - the coffeecup caustic. Note that this level set is independent of the angle $\beta$
of incidence of the incoming light.
 
The focal surface of the $k^{\mbox{th}}$ reflection of a point source also has a symmetry:

\begin{Cor}
The focal surface generated by the k$^{th}$ reflection of a point source is
invariant under translation along the cylindrical
axis.  
\end{Cor}
\begin{pf}
This follows from the fact that $x^1$ and $x^2$ of the surface in
Main Theorem 2 are independent of $u$. 
\end{pf}

This symmetry is not shared by the reflected wavefront itself - just it's focal surface. In Figure 7 we 
illustrate the 1$^{\mbox{st}}$ reflected wavefront - the lack of translational symmetry is obvious.

\vspace{0.1in}

\setlength{\epsfxsize}{4.5in}
\begin{center}
   \mbox{\epsfbox{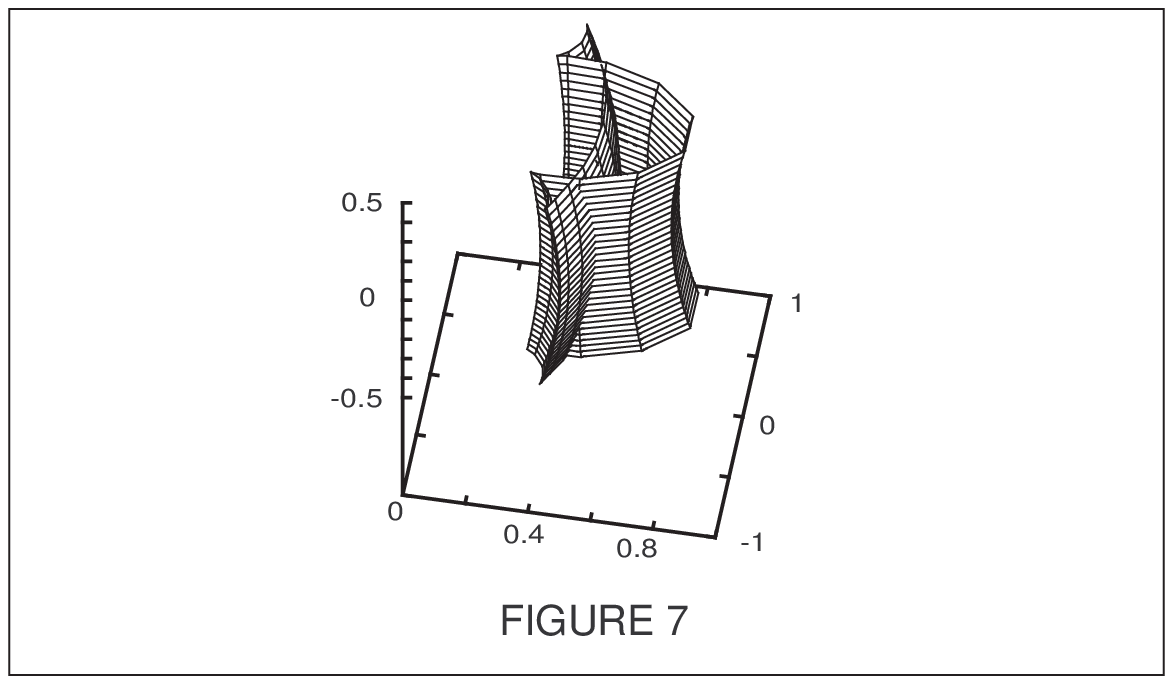}}
\end{center}

\vspace{0.1in}

In the sequence of pictures in Figure 8
we show the evolution of the level sets of the 1$^{\mbox{st}}$ 
focal set as the distance of the source decreases. The cylinder is shown by the heavy circle.

\vspace{0.1in}

\setlength{\epsfxsize}{5.0in}
\begin{center}
   \mbox{\epsfbox{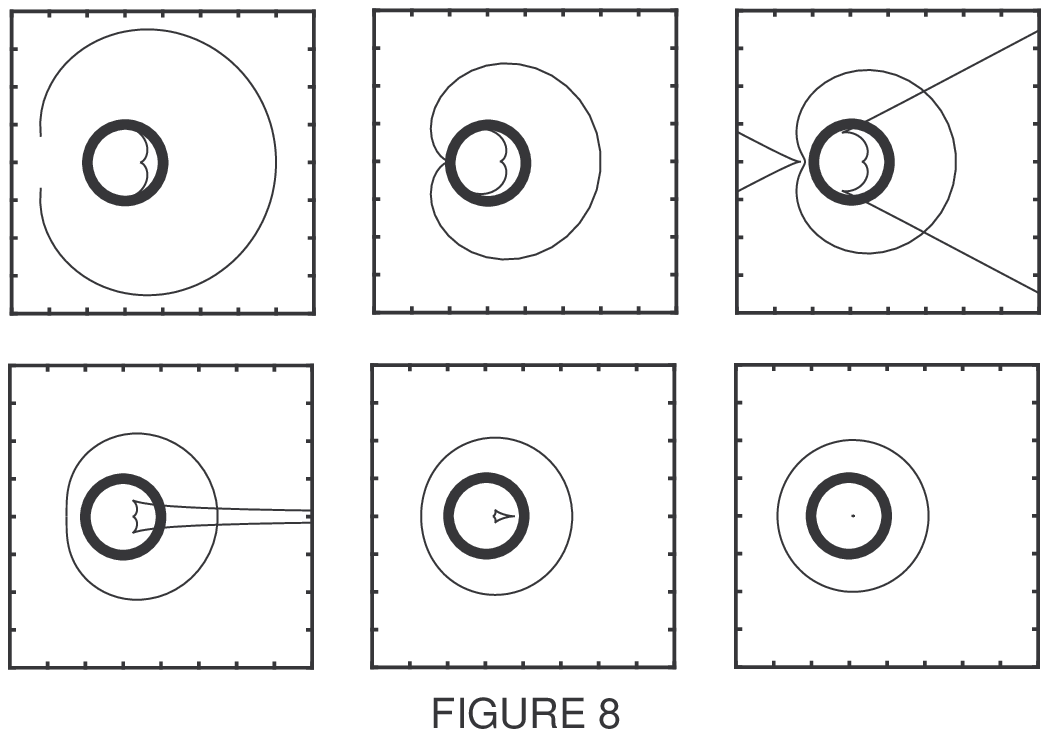}}
\end{center}

\vspace{0.1in}

The focal curve lies entirely outside of the cylinder and so is not seen in reality. 
The cross-section of the focal surface inside the cylinder for varying values
of $l/a$ is illustrated in Figure 9. The parallel wave limit,
given by Proposition \ref{p:cylpl}, is also indicated with a broken line.

\vspace{0.1in}

\setlength{\epsfxsize}{4.5in}
\begin{center}
   \mbox{\epsfbox{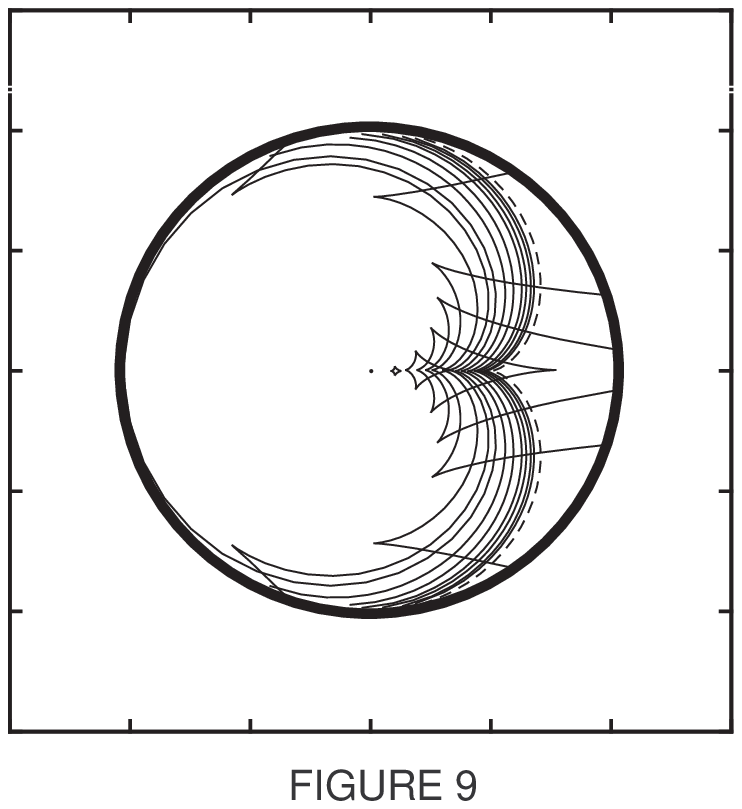}}
\end{center}

\vspace{0.1in}

In Figure 10 we compare the higher reflection caustics for
varying values of $l/a$.

\vspace{0.1in}

\setlength{\epsfxsize}{5.0in}
\begin{center}
   \mbox{\epsfbox{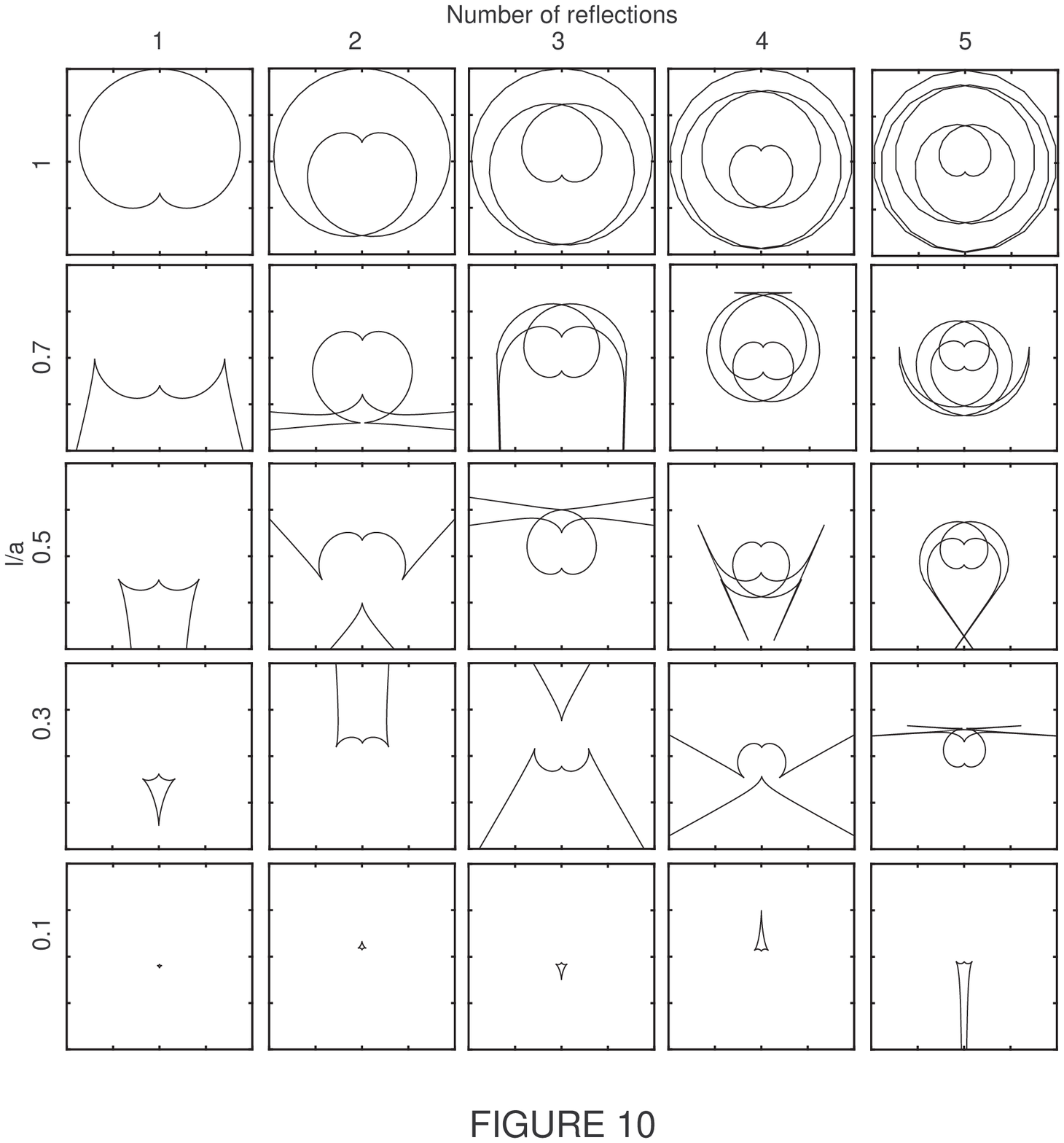}}
\end{center}

\vspace{0.1in}

These higher reflection caustics are in fact visible physically. At each reflection some of the light intensity
is lost and one expects to see a series of overlapping caustics of lessening brightness. In fact, the detailed profile of 
light intensity near a caustic varies in ways that geometric optics does not model well. Nonetheless, 
the accompanying plate is 
a photograph of the caustics formed by a 7cm diameter brass cylinder and agrees well
with the geometric optics approximation. The photograph, which was taken by the first author in collaboration with Grace 
Weir, shows the first and second reflection caustic formed by a light source at $l/a=1$ (compare with the first two curves
on the top row of Figure 10).

\end{document}